\documentclass[%
reprint,
amsmath,amssymb,
aps,
nofootinbib
]{revtex4-2}

\usepackage{graphicx}
\usepackage{dcolumn}
\usepackage{bm}
\usepackage{hyperref}
\usepackage{tikz-cd}

\usepackage{amsthm}

\newtheorem{theorem}{Theorem}
\newtheorem{lemma}{Lemma}
\newtheorem{prop}{Proposition}

\begin{document}
	\title{Optimal alignment of Lorentz orientation and generalization to matrix Lie groups}
	
	\author{Congzhou M Sha}
	\email{cms6712@psu.edu}
	\affiliation{Penn State College of Medicine, 500 University Dr, Hershey, PA, USA, 17033}
	
	\date{\today}
	
	\begin{abstract}
	There exist elegant methods of aligning point clouds in $\mathbb R^3$. Unfortunately, these methods fail to generalize to the case of Minkowski space, as we will show. Instead, we propose two solutions to the following problem: given inertial reference frames $A$ and $B$, and given (possibly noisy) measurements of a set of 4-vectors $\{v_i\}$ made in those reference frames with components $\{v_{A,i}\}$ and $\{v_{B,i}\}$, find the optimal Lorentz transformation $\Lambda$ such that $\Lambda v_{A,i}=v_{B,i}$. The first method is direct least squares optimization through a parametrization of $SO(3,1)_+$ in terms of the familiar boost and rotation vectors. The second method takes a detour through the Lorentz algebra; in addition to being conceptually simple and possessing a computational advantage over the first method, it can easily be generalized to the alignment of vector representations in other matrix Lie groups.
	\end{abstract}
	
	\maketitle
	
	\section{Introduction}
	\par The problem of finding the optimal rigid alignment between two point clouds in $\mathbb R^3$ is well-known, and is also referred to as the \textbf{orthogonal Procrustes problem}. Two particularly elegant and efficient solutions exist: the Kabsch algorithm \cite{Kabsch1976} and the Horn algorithm \cite{Horn1987}. These algorithms require as input two lists of vector components $\{v_{A,i}\}$ and $\{v_{B,i}\}$, which share a common origin, where $v_{A,i}$ is to be mapped to $v_{B,i}$. The Kabsch algorithm easily extends to arbitrary dimensions $N$, and produces an element of $SO(N)$. In contrast, the Horn algorithm relies on the geometric connection between unit quaternions and 3D rotations to rephrase the problem into a principal eigenvalue/eigenvector problem, and so it can only be used for $SO(3)$.
	\par Both methods rely on the positive definiteness of the Euclidean metric, or equivalently that $SO(N)$ is compact, ensuring that the resulting optimization problem is convex. In contrast, the Minkowski inner product (we use the convention $\eta_{\mu\nu}=\text{diag}(-,+,+,+)$ in this work) induces a norm which is indefinite, and the proper orthochronous Lorentz group $SO(3,1)_+$ is noncompact due to the Lorentz boosts.
	\par In the appendices, we show that neither of the algorithms above can be applied directly to the problem of alignment in $SO(3,1)_+$, and an alternative approach is needed to find the optimal Lorentz transformation from one set of 4-vectors to another.
	\par In this work, we present two methods of performing this alignment. The first method is traditional least squares, constrained optimization using Newton-Raphson iteration, which has the advantage of being highly numerically robust, but is computationally intensive. The second method is by solving for an least squares linear transformation with the Moore-Penrose pseudoinverse (which does not guarantee a Lorentz matrix such that $L\eta L^T=\eta$), and projecting that linear transformation onto the nearest element of the Lorentz group via a detour through the Lie algebra $\mathfrak{so}(3,1)$. To our knowledge, this problem has not been addressed previously in the literature.
	\par Note that in flat space, points and tangent vectors are conflated. In this work, we assume that all inputs are true tangent vectors, and there is no ``centering" process required.

	\section{Proposed solutions to the Lorentz alignment problem}
	\subsection{Review of the Lorentz algebra and group}
	To review, matrix Lie groups such as $SO(3,1)_+$ are related to a Lie algebra, e.g. $\mathfrak{so}(3,1)$, via the matrix logarithm and matrix exponential functions. The complex matrix logarithm exists if the matrix is invertible (Theorem 2.10 in \cite{Hall2013a}). The Lorentz algebra is the set of matrices:
	\begin{align}
		\label{eq:lorentzalg}
		A=\begin{bmatrix}
			0&\zeta^1&\zeta^2&\zeta^3\\
			\zeta^1&0&-\theta^3&\theta^2\\
			\zeta^2&\theta^3&0&-\theta^1\\
			\zeta^3&-\theta^2&\theta^1&0
		\end{bmatrix}\in \mathfrak{so}(3,1),
	\end{align}
	where $\boldsymbol{\zeta}=[\zeta^1,\zeta^2,\zeta^3]$ is the boost vector and $\boldsymbol\theta=[\theta^1,\theta^2,\theta^3]$ is the rotation vector, and both are real. The matrix exponential of this set of matrices (i.e. the Lorentz analogue of Rodrigues' rotation formula \cite{Dai2015}) has recently been explicitly computed \cite{Haber2024}:
	\begin{align}
		\label{eq:matrixexpLorentz}
		\Lambda(\boldsymbol\zeta,\boldsymbol\theta)=\exp A=\frac{1}{a^2+b^2}\sum_{i=0}^3 f_i(a,b)A^i,\\
		a^2=\frac12\left[|\boldsymbol\theta|^2-|\boldsymbol\zeta|^2+\sqrt{\left(|\boldsymbol\theta|^2-|\boldsymbol\zeta|^2\right)^2+4\left(\boldsymbol\theta\cdot\boldsymbol\zeta\right)^2}\right],\\
		b^2=\frac12\left[|\boldsymbol\zeta|^2-|\boldsymbol\theta|^2+\sqrt{\left(|\boldsymbol\theta|^2-|\boldsymbol\zeta|^2\right)^2+4\left(\boldsymbol\theta\cdot\boldsymbol\zeta\right)^2}\right],\\
		\label{eq:haberf1}
		f_0(a,b)=b^2\cos a+a^2\cosh b,\\
		f_1(a,b)=\frac{b^2}a\sin a+\frac{a^2}b\sinh b,\\
		f_2(a,b)=\cosh b-\cos a,\\
		\label{eq:haberf4}
		f_3(a,b)=\frac{\sinh b}a-\frac{\sin a}a.
	\end{align}
	Haber notes that Eqs (\ref{eq:haberf1})-(\ref{eq:haberf4}) are independent of the overall signs of $a$ and $b$, which is evident by noting that $\cos$ and $\cosh$ are even, whereas $\sin$ and $\sinh$ are odd. Since it is a continuous function, the matrix exponential reproduces only the proper orthochronous part of $SO(3,1)$, which is what we would presumably desire of an alignment algorithm.
	\subsection{Method 1: Direct nonlinear optimization}
	For naive least squares optimization of alignment, we write down the objective function, using the Frobenius norm $\|\cdot\|$ to match numerical components:
	\begin{align}
		f(\boldsymbol\zeta,\boldsymbol\theta)=\sum_i\|v_{B,i}-\Lambda(\boldsymbol\zeta,\boldsymbol\theta)v_{A,i}\|^2,
	\end{align}
	which can be minimized either through gradient descent or another method \cite{Burden2016}, using SciPy's general solvers for this task \cite{Virtanen2020a}. We will refer to this method as the direct minimization method. Note that iterative solvers require a suitable initial value, for which we used the identity matrix.
	\subsection{Method 2: Projection onto the Lorentz group after unconstrained solving}
	A common method to solve problems of the form 
	\begin{equation}
		Ax=b
	\end{equation}
	in a least squares sense is to use the Moore-Penrose inverse (or pseudoinverse) $A^+$:
	\begin{equation}
		x=A^+b,
	\end{equation}
	where $A^+=(A^TA)^{-1}A^T$ when $(A^TA)^{-1}$ exists \cite{Burden2016}. Consider now the $4\times n$ matrices $X_{ji}=(v_{A,i})_j$ and $Y_{ji}=(v_{B,i})_j$ (rows index the dimension, and columns index the vectors to be aligned), and the equation
	\begin{align}
		\label{eq:init}
		\Lambda_0 X=Y.
	\end{align}
	The Moore-Penrose inverse gives a least squares estimate for a linear transformation $\Lambda_0$ which transforms the vectors $v_{A,i}$ to $v_{B,i}$, but which is not necessarily a Lorentz transform. In other words, we assume that there is some measurement error in $X$ and $Y$ which prevents $\Lambda_0$ from satisfying $\Lambda_0^T\eta\Lambda_0=\eta$. To project $\Lambda_0$ onto a nearby true Lorentz transformation, we take a matrix logarithm:
	\begin{equation}
		\label{eq:matlog}
		l_0=\log \Lambda_0=\begin{bmatrix}
			m_{11}&m_{12}&m_{13}&m_{14}\\
			m_{21}&m_{22}&m_{23}&m_{24}\\
			m_{31}&m_{32}&m_{33}&m_{34}\\
			m_{41}&m_{42}&m_{43}&m_{44}
		\end{bmatrix},
	\end{equation}
	and consider an objective function of the Frobenius norm
	\begin{align}
		\label{eq:obj2}
		f(l)=\|l-l_0\|^2,\\
		\label{eq:obj3}
		l\in\mathfrak{so}(3,1).
	\end{align}
	In essence, this objective function forces $l_0$ to be of the form given in Eq (\ref{eq:lorentzalg}). Specifically, (1) we set the diagonal of $l_0$ to be $0$, (2) symmetrize the first row and first column, and (3) antisymmetrize the bottom right $3\times3$ block:
	\begin{equation}
		\label{eq:l_alg}
		l=\frac12\begin{bmatrix}
			0&m_{12}+m_{21}&m_{13}+m_{31}&m_{14}+m_{41}\\
			m_{12}+m_{21} & 0 & m_{23}-m_{32} &m_{24}-m_{42}\\
			m_{13}+m_{31} & m_{32}-m_{23} & 0 & m_{34}-m_{43}\\
			m_{14}+m_{41} & m_{42}-m_{24} & m_{43}-m_{34} & 0
		\end{bmatrix}
	\end{equation}
	This solution is seen since minimization of Eq (\ref{eq:obj2}) is equivalent to componentwise minimization. Since $l\in \mathfrak{so}(3,1)$, the diagonal components of $l$ must be $0$, and by symmetry (and anti-symmetry), the off-diagonal elements minimize Eq (\ref{eq:obj2}) at the averages of $\boldsymbol\theta_i$ and $\boldsymbol\zeta_i$. Finally, we take $\Lambda=\exp l$ using Eqs (\ref{eq:matrixexpLorentz})-(\ref{eq:haberf4}) as the optimal Lorentz transformation.
	\par We have thus defined a projection operator
	\begin{equation}
		\text{proj}_{\mathfrak{so}(3,1)}(l_0)=l,
	\end{equation}
	which is a linear operator since each entry of $l$ is a linear combination of entries of $l_0$.
	\begin{lemma}
		\label{lemma:1}
		To linear order, minimization in the Lie algebra in Eq (\ref{eq:obj2}) and (\ref{eq:obj3}) corresponds to minimization of the following function of the Frobenius norm
		\begin{align}
			\label{eq:obj4}
			f(\Lambda)=\|\log \left(\Lambda^{-1}\Lambda_0\right)\|^2,\\
			\Lambda\in SO(3,1).
		\end{align}
	\end{lemma}
	\begin{proof}
		We recall the Baker-Campbell-Hausdorff formula (\cite{Hall2013a}):
		\begin{equation}
			\log\left[\exp X \exp Y\right]=X+Y+\frac12[X,Y]+\cdots,
		\end{equation}
		where $[X,Y]=XY-YX$. Since $\Lambda=\exp(l)$ and $\Lambda_0=\exp(l_0)$,
		\begin{align}
			\log\left(\Lambda^{-1}\Lambda_0\right)=\log\left(\exp(-l)\exp(l_0)\right)\\\approx -l+l_0+\frac12[-l,l_0]
		\end{align}
	\end{proof}
	\begin{theorem}
		\label{thm:minequiv}
		Let $\Lambda_{\text{alg}}=\exp(l)$ be the Lorentz transformation obtained via the Lie algebra method, where $l=\text{proj}_{\mathfrak{so}(3,1)}(\log\Lambda_0)$ is the orthogonal projection (with respect to the Frobenius inner product) of the matrix logarithm onto the Lie algebra $\mathfrak{so}(3,1)$. Let $\Lambda_{\text{GT}}\in SO(3,1)_+$ be the ground truth Lorentz transformation. Under sufficiently small measurement noise ($\|\Delta\|\ll\|\Lambda_{\text{GT}}\|$ as defined below), the Lie algebra method recovers the ground truth with error
		\begin{equation}
			\|\Lambda_{\text{alg}}-\Lambda_{\text{GT}}\|=O\left(\|\Lambda_{\text{GT}}^{-1}EX^+\|\left(1+\|\log\Lambda_{\text{GT}}\|\right)\right).
		\end{equation}
	\end{theorem}
	\begin{proof}
		Let $\Lambda_{\text{GT}}\in SO(3,1)_+$ be the ground truth Lorentz transformation, and let there be noise $E$ in the measurement of $Y$:
		\begin{equation}
			Y=\Lambda_{\text{GT}}X+E.
		\end{equation}
		Using Eq (\ref{eq:init}),
		\begin{equation}
			Y=\Lambda_0X=\Lambda_{GT}X+E,
		\end{equation}
		Assuming that the Moore-Penrose inverse $X^+$ exists,
		\begin{equation}
			\Delta=\Lambda_0-\Lambda_{\text{GT}}=EX^+.
		\end{equation}
		If $\|\Delta\|\ll \|\Lambda_{\text{GT}}\|$ for some matrix norm $\|\cdot\|$, then we may write
		\begin{align}
			l_0 &= \log(\Lambda_{\text{GT}}+\Delta)\nonumber\\
			&= \log\left[\Lambda_{\text{GT}}\left(I+\Lambda_{\text{GT}}^{-1}\Delta\right)\right]\nonumber\\
			&\approx \log\Lambda_{\text{GT}}+\log\left(I+\Lambda_{\text{GT}}^{-1}\Delta\right)\nonumber\\
			&\approx \log\Lambda_{\text{GT}}+\Lambda_{\text{GT}}^{-1}\Delta-\frac12\left(\Lambda_{\text{GT}}^{-1}\Delta\right)^2+\cdots
		\end{align}
		Let $l_{\text{GT}}=\log\Lambda_{\text{GT}}\in\mathfrak{so}(3,1)$. The Lie algebra method projects $l_0$ onto $\mathfrak{so}(3,1)$, producing $l$, which when exponentiated gives $\Lambda_{\text{alg}}=\exp(l)$. Since $l_{\text{GT}}\in\mathfrak{so}(3,1)$, we have
		\begin{align}
			l=\text{proj}_{\mathfrak{so}(3,1)}(l_0)\nonumber\\
			=\text{proj}_{\mathfrak{so}(3,1)}(l_{\text{GT}})+\text{proj}_{\mathfrak{so}(3,1)}\left(\Lambda_{\text{GT}}^{-1}\Delta\right) + O(\|\Delta\|^2)\nonumber\\
			=l_{\text{GT}} + \text{proj}_{\mathfrak{so}(3,1)}\left(\Lambda_{\text{GT}}^{-1}\Delta\right) + O(\|\Delta\|^2),
		\end{align}
		where in the second line we have used the linearity of the projection operator.
		\\\\
		Let $\delta l = \text{proj}_{\mathfrak{so}(3,1)}(\Lambda_{\text{GT}}^{-1}\Delta)$ denote the projected perturbation. By Lemma \ref{lemma:1}, the Baker-Campbell-Hausdorff formula gives
		\begin{align}
			\Lambda_{\text{alg}} &= \exp(l) = \exp(l_{\text{GT}} + \delta l + O(\|\Delta\|^2))\nonumber\\
			&\approx \exp(l_{\text{GT}})\exp(\delta l)\exp\left(\frac{1}{2}[l_{\text{GT}}, \delta l]\right).
		\end{align}
		Since the projection is bounded by $\|\delta l\|=\|\text{proj}_{\mathfrak{so}(3,1)}(\Lambda_{\text{GT}}^{-1}\Delta)\|\leq\|\Lambda_{\text{GT}}^{-1}\Delta\|=O(\|EX^+\|)$, and using $\exp(X)=I+X+O(\|X\|^2)$, we have
		\begin{equation}
			\exp(\delta l) = I + O(\|\Lambda_{\text{GT}}^{-1}EX^+\|).
		\end{equation}
		The commutator term can be bounded using submultiplicativity of the matrix norm:
		\begin{align}
			\|[l_{\text{GT}}, \delta l]\| &\leq 2\|l_{\text{GT}}\|\cdot\|\delta l\|\nonumber\\
			&\leq 2\|\log\Lambda_{\text{GT}}\|\cdot\|\Lambda_{\text{GT}}^{-1}EX^+\|.
		\end{align}
		Combining with the leading-order term $\|\delta l\|\leq\|\Lambda_{\text{GT}}^{-1}EX^+\|$, we obtain the stated error bound:
		\begin{align}
			\|\Lambda_{\text{alg}}-\Lambda_{\text{GT}}\| = O\left(\|\Lambda_{\text{GT}}^{-1}EX^+\|(1+\|\log\Lambda_{\text{GT}}\|)\right).
		\end{align}
	\end{proof}
	To refine the error in Theorem \ref{thm:minequiv}, we recall that the Moore-Penrose inverse can be constructed from its singular value decomposition \cite{Strang2003}, and that the spectral radius of $X^+$ therefore is $\frac{1}{\sigma_{\text{min}}(X)}$, where $\sigma_{\text{min}}$ is the lowest magnitude singular value of $X$. The factor of $X^+$ in the error bound therefore contributes a factor of $\frac{1}{\sigma_{\text{min}}(X)}$ to the overall error. The factor $(1+\|\log\Lambda_{\text{GT}}\|)$ indicates that the error grows when the ground truth transformation is far from the identity (i.e., involves large boosts or rotations).
	\par If $X$ is rank deficient, one of its singular values must be zero, which may prove a barrier in computing the Moore-Penrose inverse. In the case of minimizing Eq (\ref{eq:cost_lorentz}), this indicates that not enough independent measurements were performed, and there may be infinitely many solutions via the direct minimization method. Such an issue can be overcome either by making more measurements or by offsetting the zero-valued singular values by a small real number $\epsilon$. The Lie algebra method may encounter difficulties in Eq (\ref{eq:matlog}).
	\par While the matrix logarithm always exists over $\mathbb{C}$ for any invertible matrix (since $\mathbb{C}$ is algebraically closed), and $\exp(\log \Lambda_0)=\Lambda_0$ always holds by definition, the logarithm is not unique due to branch cuts. However, standard numerical implementations (such as SciPy's \verb|scipy.linalg.logm|) reliably compute the principal branch, which yields a real matrix for matrices in or near $SO(3,1)_+$. Additionally, the choice of branch cut does not affect the final result, since we only use the matrix exponential of any matrix logarithms.
	\par When the noise in $\Lambda_0$ is sufficiently large that it lies far from $SO(3,1)_+$, the computed logarithm may not be close to any element of $\mathfrak{so}(3,1)$, and the projection step may produce spurious results. We show empirically in the Results that such pathologies were not encountered in what one might consider typical usage of this method.
	\par This method of projection extends easily to all other matrix Lie groups: one simply proposes and minimizes an objective function on the Lie algebra, taking care that the appropriate singular values are nonzero, and checking that measurement noise is within tolerable bounds. In certain cases the optimal solution in the Lie algebra can be written down explicitly, as in Eq (\ref{eq:l_alg}). Extra consideration is needed if the Lie group is not simply connected and one desires elements other than from the connected component containing the identity. We refer to this method as the Lie algebra method.
	\section{Results}
	We implemented both methods in Python, benchmarking accuracy and timing. We used NumPy \cite{Harris2020} for array operations, and SciPy \cite{Virtanen2020a} matrix logarithms, matrix exponentials, and nonlinear solving. All testing was performed on an M3 MacBook Air. 
	\subsection{Ideal conditions with no noise}
	We began with a set of four time-like 4-vectors:
	\begin{align}
		X=\begin{bmatrix}
			1&0&0&0\\
			\sqrt2&1&0&0\\
			\sqrt2&0&1&0\\
			\sqrt2&0&0&1
		\end{bmatrix}
	\end{align}
	and a simple boost along $x$ of the form
	\begin{align}
		\Lambda=\begin{bmatrix}
			\gamma&-\beta\gamma&0&0\\
			-\beta\gamma&\gamma&0&0\\
			0&0&1&0\\
			0&0&0&1
		\end{bmatrix},
	\end{align}
	with $\beta=0.3$ and $\gamma=\frac{1}{\sqrt{1-\beta^2}}$, and verified that both methods could reproduce $\Lambda$, given $A$ and $B=A\Lambda^T$. This simple sanity check showed that both methods indeed work. The use of the explicit form of the Lorentz boost given by Haber in Eqs (\ref{eq:matrixexpLorentz})-(\ref{eq:haberf4}) ensured that the resulting matrix indeed had determinant close to $1$. However, the Lie algebra method was slightly more effective than direct least squares optimization via the default solver in \verb|scipy.optimize.minimize| in ensuring that the components of $A\Lambda_{\text{est}}^T$ were close to $A\Lambda^T$, as demonstrated in the attached Code and Data. For the Lie algebra method, we tested brute force matrix exponentiation via \verb|scipy.linalg.expm| compared to Haber's explicit computation, and found that both computations resulted in $\det \Lambda\approx 1$ to 14 decimal places.
	\par The timing between the two Lie algebra methods was also extremely similar, with each iteration taking 700-800 $\mu$s, with standard deviation of 20-30 $\mu$s. In contrast, the default solver for \verb|scipy.optimize.minimize| took an average of 24 ms, which was $\sim30\times$ as slow as the Lie algebra method.
	\subsection{Varying the number of vectors given, noise, and averaging over arbitrary Lorentz transformations}
	In Figs \ref{fig:l2} and \ref{fig:linf}, we compare the absolute Frobenius ($L^2$) and max ($L^\infty$) norms of the estimated Lorentz transformation matrix ($\Lambda_{\text{est}}$) against the true Lorentz transformation ($\Lambda_{\text{GT}}$). We generated $1,000$ random Lorentz transformations using
	\begin{align}
		\zeta^i\sim\mathcal{N}(0,0.2),\\
		\theta^i\sim\mathcal{N}(0, 1).
	\end{align}
	\begin{widetext}
	\begin{figure*}
		\includegraphics[width=\linewidth]{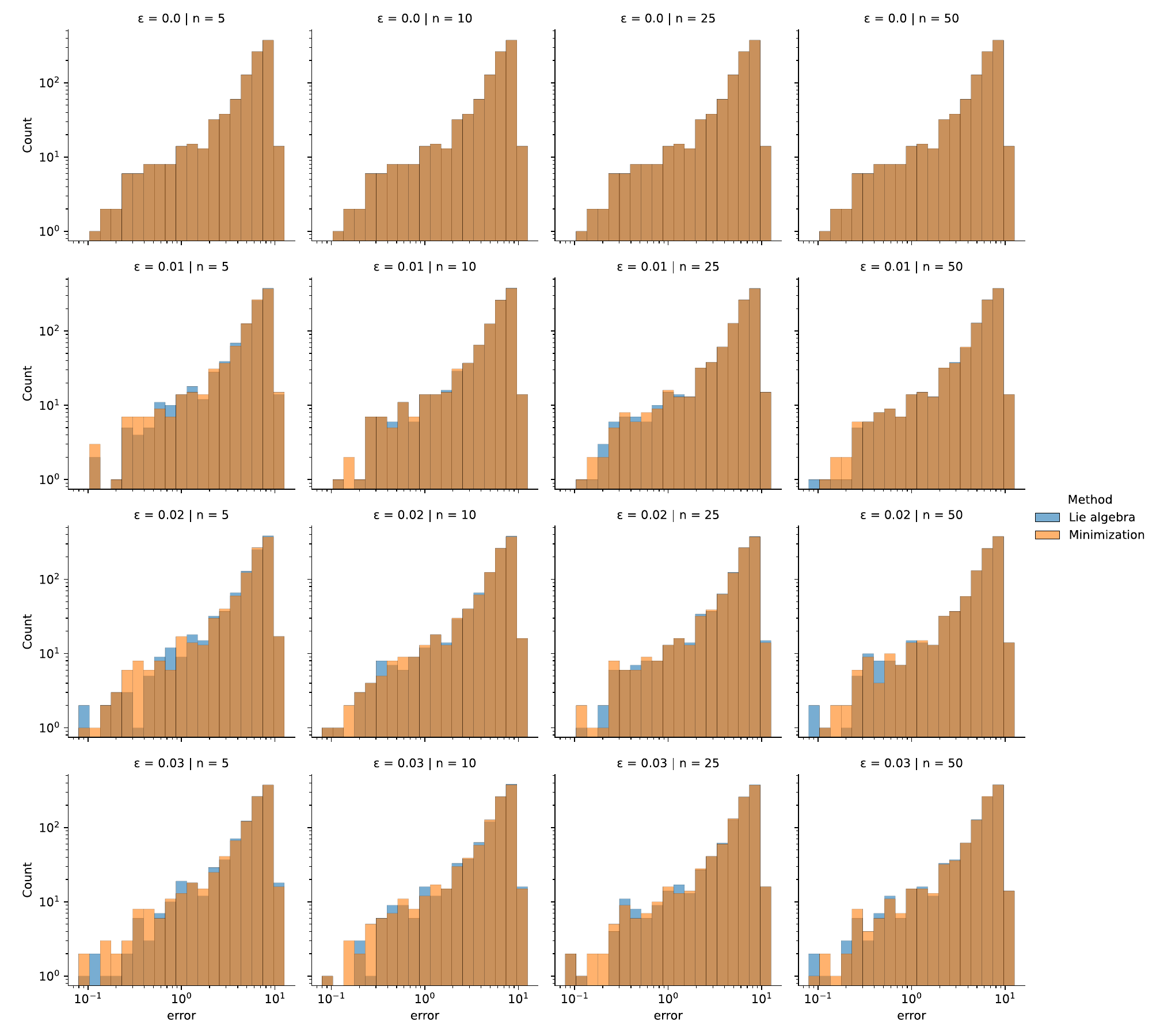}
		\caption{The absolute Frobenius ($L^2$) norm of the estimated Lorentz matrix from frame $A$ to frame $B$ as compared to the actual Lorentz matrix ($\|\Lambda_{\text{est}}-\Lambda_{\text{GT}}\|_2$). The rows represent varying the noise $\epsilon$ present in the measurement of the vectors in frame $B$. The columns represent varying the number $n$ of vectors upon which the estimation was performed. Evidently, the Lie algebra method and the direct minimization method are equivalent in accuracy across a range of testing conditions, since their absolute error distributions largely overlap.}
		\label{fig:l2}
	\end{figure*}
		\begin{figure*}
	\includegraphics[width=\linewidth]{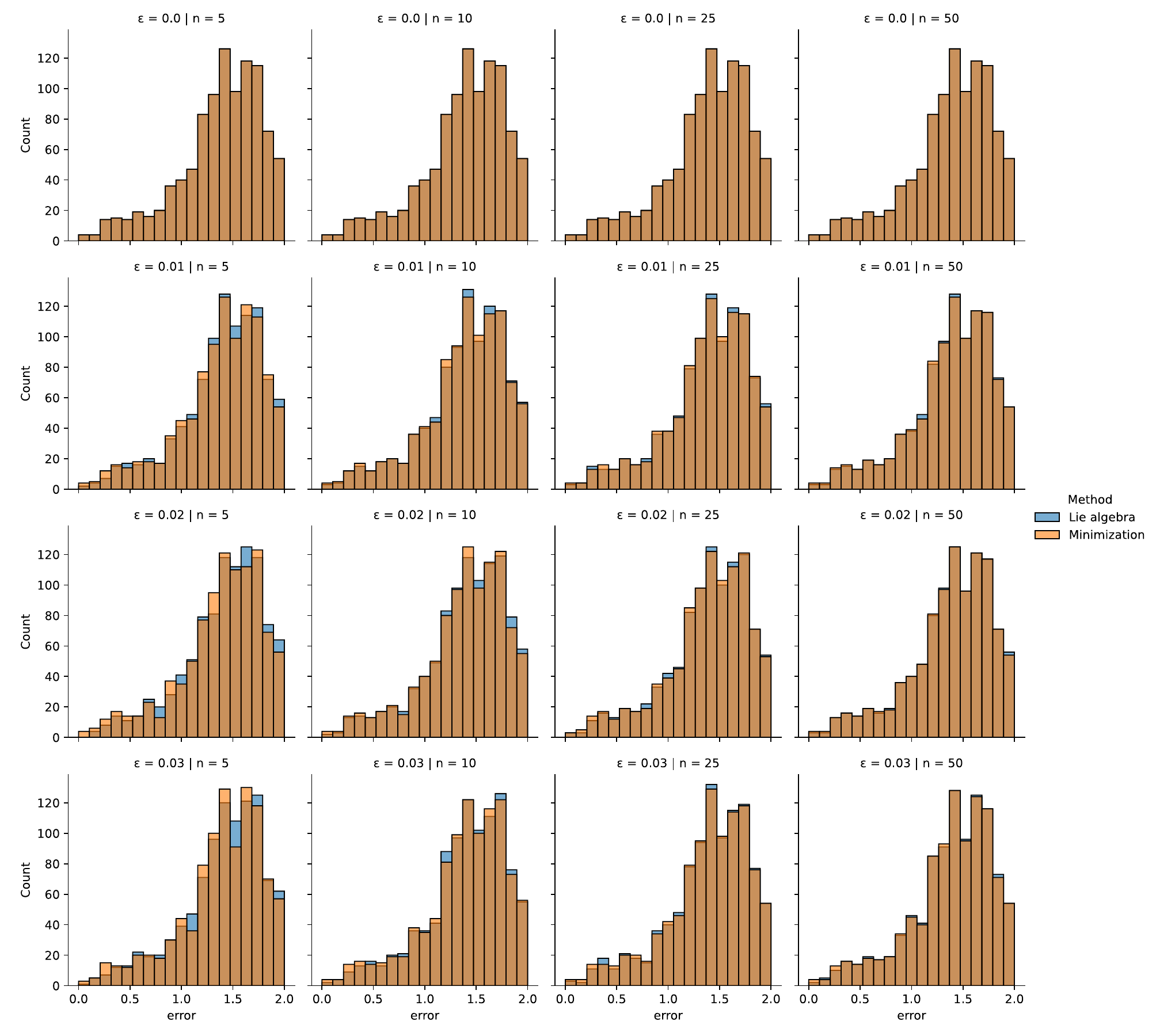}
	\caption{The absolute max ($L^\infty$) norm of the estimated Lorentz matrix from frame $A$ to frame $B$ as compared to the actual Lorentz matrix ($\|\Lambda_{\text{est}}-\Lambda_{\text{GT}}\|_text{max}$). The rows and columns represent the same variables as in Fig \ref{fig:l2}.}
	\label{fig:linf}
\end{figure*}
	\end{widetext}
	Note that the joint distribution for $\boldsymbol\zeta$ and $\boldsymbol\theta$ is a spherically symmetric product of normal distributions. We generated $n$ random 4-vectors with Minkowski norm $-1$ by sampling $x,y,z\sim\mathcal{N}(0,0.3)$ and setting $t=\sqrt{1+x^2+y^2+z^2}$, and used these vectors as our $v_{A,i}$. We introduced noise in the vector components as measured in frame $B$ by adding random noise $\Delta x,\Delta y,\Delta z\sim\mathcal{N}(0,\epsilon)$ to $x,y,z$, normalizing the time component for $x+\Delta x,y+\Delta y,z+\Delta z$ again, and then boosting by $\Lambda$ to these vectors.
	\par We observe that the overall distribution of both the absolute maximum element-wise error and the overall $L^2$ error is quite similar between the two proposed methods (Figs \ref{fig:l2}, \ref{fig:linf}).
	
	\section{Discussion}
	In this work, we have proposed two methods of solving the Lorentz alignment problem. In the appendices, we showed that the principles behind the Kabsch algorithm and Horn algorithm for rigid rotation alignment in Euclidean $\mathbb R^3$ do not easily extend to the Lorentz geometry. A Kabsch-like approach is not feasible unless a Lorentz singular value decomposition is available (Proposition \ref{thm:kabsch}), the development of which is beyond the scope of this work. The Horn-like approach is infeasible because of the unbounded nature of the unit biquaternions representation of the Lorentz group, as compared to the bounded real quaternion representation of $SO(3)$; this obstacle prevents the optimization problem from linearizing into a principal eigenvalue/eigenvector problem (Proposition \ref{thm:horn}).
	\par Instead, we can either perform constrained optimization using gradient descent or related methods (the direct minimization method) or project the least-squares linear operator via the pseudoinverse onto the Lorentz group via a second least squares minimization in the Lie algebra (the Lie algebra method). We have shown that at least in Python, the second approach is more performant, and better reproduces the answer in the case that an exact solution exists. We did not vary the direct minimization method used by SciPy's default solver, though ideally an optimal solution would not be dependent on such as choice. The explicit formula relating the boost and rotation vectors to the corresponding Lorentz transformation was necessary in the direct minimization method, whereas it was optional in the Lie algebra method with little performance difference between it and brute force matrix exponentiation.
	\par The Lie algebra method presented here technically also solves the orthogonal Procrustes problem by restricting our attention to the bottom right $3\times 3$ matrices in the formulae above. However, it extends to arbitrary matrix Lie groups, unlike the Kabsch algorithm which only works for $O(N)$, and the Horn algorithm which relies on a quaternion representation with bounded components. Therefore we feel that this work is valuable in the overall scheme of manifold alignment. An application of Lorentz alignment is finding the connection between free-falling reference frames in smooth lattice general relativity \cite{Brewin2011}.
	\par The direct minimization method requires repeated evaluation of the matrix exponential in Eqs (\ref{eq:matrixexpLorentz})-(\ref{eq:haberf4}), on the order of tens to hundreds of times, in order to perform gradient descent or some other method of iterative optimization. In contrast, the Lie algebra method requires a single evaluation of the matrix logarithm which may be computed iteratively \cite{Al-Mohy2012}, and a single evaluation of the matrix exponential. Depending on the convergence speed of the matrix logarithm in comparison to the matrix exponential, the Lie algebra method likely outperforms the direct minimization method for other matrix Lie groups, though further numerical experiments are necessary to verify this hypothesis. In both methods, we used an explicit formula to compute the matrix exponential, however this is also not essential since the matrix exponential can be computed through eigendecomposition. Again, further numerical experiments would be required to compare the performance of explicit formulae for the matrix exponential versus exponentiation via eigendecomposition for other matrix Lie groups. Finally, implementation of these algorithms in other languages such as Julia and C would enable more comprehensive benchmarks.
	\section*{Code and data availability}
	All code needed to reproduce the data and figures in this work are available at Zenodo (doi: \href{https://zenodo.org/records/15724187}{10.5281/zenodo.15724187}).
	\section*{Acknowledgments}
	The author would like to thank the Penn State College of Medicine's Medical Scientist Training Program for its support of his career.
	\section*{Conflict of interest statement}
	The author reports no potential conflicts of interest.
	\appendix
	
	\section{Quaternions}
		\label{sec:quaternions}
	Quaternions (often represented by the symbol $\mathbb H$) can be thought of as a four-dimensional real vector space with basis $\{\mathbf 1,\mathbf i,\mathbf j,\mathbf k\}$, with quaternion addition defined by componentwise addition. We will represent quaternions as $q=(q_0,q_1,q_2,q_3)$, where $q_0,q_1,q_2,q_3\in\mathbb R$. We denote the imaginary unit in $\mathbb C$ by $h^2=-1$. 	
	The operation of quaternion conjugation is defined by
	\begin{equation}
		\label{eq:quatconj}
		q^*=(q_0,-q_1,-q_2,-q_3),
	\end{equation}
	Further algebraic properties define the operation of quaternion multiplication:
	\begin{eqnarray}
		\mathbf i\mathbf j=\mathbf k,\\
		\mathbf j\mathbf k=\mathbf i,\\
		\mathbf k\mathbf i=\mathbf j,\\
		\mathbf i^2=\mathbf j^2=\mathbf k^2=-\mathbf 1,
	\end{eqnarray}
	so that
	\begin{align}
		(p_0,p_1,p_2,p_3)&(q_0,q_1,q_2,q_3)\nonumber\\=(
		& p_0q_0-p_1q_1-p_2q_2-p_3q_3,\nonumber\\
		& p_0q_1 + p_1q_0 + p_2q_3 - p_3q_2,\nonumber\\
		& p_0q_2 - p_1q_3 + p_2q_0 + p_3q_1,\nonumber\\
		& p_0q_3 + p_1q_2 - p_2q_1 + p_3q_0).
		\label{eq:quatmul}
	\end{align}
	Considering $q$ as a column vector, quaternion multiplication $pq$ can also be described as left multiplication of that column vector by a matrix
	\begin{equation}
		\label{eq:quatmatmul}
		pq=\begin{bmatrix}
			p_0&-p_1&-p_2&-p_3\\
			p_1&p_0&p_3&-p_2\\
			p_2&-p_3&p_0&p_1\\
			p_3&p_2&-p_1&p_0
		\end{bmatrix}\begin{bmatrix}q_0\\q_1\\q_2\\q_3\end{bmatrix}.
	\end{equation}
	This interpretation of quaternion multiplication as matrix multiplication is used extensively by Horn. Quaternion multiplication is associative but not commutative.
	\par The quaternion norm is defined as $||q||_{\text{quat}}=qq^*=q_0^2+q_1^2+q_2^2+q_3^2$, and unit quaternions by $||q||_{\text{quat}}=1$. It is unit quaternions which are isomorphic $SU(2)$, which is a double cover of $SO(3)$. To represent rotations, a vector $\mathbf{r}\in \mathbb R^3$ is encoded as a quaternion $r=(0,r_1, r_2, r_3)$. We define the operator $U$ which undoes this encoding by $U(r)=\mathbf r$. Then the rotation of $\mathbf{r}\in \mathbb R^3$ by a matrix $R\in SO(3)$ can alternatively be performed in quaternion space by:
	\begin{align}
		\label{eq:quaternify}
		R\mathbf{r}=U(qrq^*),
	\end{align}
	where $q$ is the quaternion representing $R$, the left-hand side $R\mathbf{R}$ is a matrix multiplication, and the right-hand side $qrq^*$ is a product of three quaternions, computed using Eq (\ref{eq:quatmul}).
	\par The components of $q$ which correspond to a rotation around an axis $\mathbf u$ counterclockwise by an angle $\theta$ are:
	\begin{equation}
		\label{eq:axisanglequat}
		q=\left(\cos\frac\theta2,u_1\sin\frac\theta2,u_2\sin\frac\theta2,u_3\sin\frac\theta2\right).
	\end{equation}
	Using Eq (\ref{eq:axisanglequat}) along with Rodrigues' rotation formula \cite{Dai2015}, Eq (\ref{eq:quaternify}) can be directly verified.
		\section{Insights into pose alignment from the Kabsch and Horn algorithms}
	In this section, we review the proofs of the Kabsch and Horn algorithms, and justify why neither method is directly applicable to alignment in the Lorentz group.
	\subsection{The Kabsch and Horn algorithms}
	\begin{theorem}
		\label{thm:kabsch}
		(Kabsch algorithm)
		Let there be two lists $\{\mathbf{v}_{A,i}\}$ and $\{\mathbf{v}_{B,i}\}$ of $n$ vectors in $N$ dimensions. Let $A$ and $B$ be $N\times n$ matrices where $A_{ji}=(\mathbf{v}_{A,i})_j$ and $B_{ji}=(\mathbf{v}_{B,i})_j$. Let the objective function be
		\begin{equation}
			\label{eq:cost1}
			L(R)=\sum_i \|\mathbf{v}_{B,i}-R\mathbf{v}_{A,i}\|^2=\|B-RA\|^2,
		\end{equation}
		where $\|\cdot\|$ is the Frobenius norm. Given the singular value decomposition of $BA^T=U\Sigma V^T$. If $|UV^T|=1$, the rotation matrix $R$ which minimizes $L(R)$ is \begin{equation}R=UV^T,\end{equation}
		and otherwise \begin{equation}R=U\begin{bmatrix}1&0&0\\0&1&0\\0&0&-1\end{bmatrix}V^T.\end{equation}
	\end{theorem}
	\begin{proof}
		First, we simplify the objective function:
		\begin{align}
			L(R)=\text{tr}\left((B-RA)^T(B-RA)\right)\\
			\label{eq:cost3}
			=\text{tr}(B^TB)-2\text{tr}(A^TR^TB)+\text{tr}(A^TA).
		\end{align}
		In Eq (\ref{eq:cost1}), the left side norm is the least squares (i.e. Frobenius) norm for vectors, and the norm on the right side is the least squares norm for matrices. In Eq (\ref{eq:cost3}), the first and final terms are independent of $R$, we may drop them and seek to minimize the middle term, or equivalently after using the cyclicity property of the trace, maximize:
		\begin{align}
			\label{eq:cost4}
			f(R)=\text{tr}(R^TBA^T).
		\end{align}
		The key insight in the Kabsch algorithm is to perform singular value decomposition on the matrix $H=BA^T=U\Sigma V^T$, where $U, V\in O(N)$, and $\Sigma=\text{diag}(\sigma_1, \sigma_2,\cdots)$ contains the singular values.
		\begin{align}
			f(R)=\text{tr}(R^TH)\\
			=\text{tr}(R^TU\Sigma V^T)\\
			\label{eq:simp1}
			=\text{tr}(V^TR^TU\Sigma).
		\end{align}
		Since $V^T,R^T,U\in O(N)$, the resulting product is $M=V^TR^TU\in O(N)$. Furthermore, $|M_{ii}|\leq 1$ and $\text{tr}(M_{ii})\leq N$.
		\begin{align}
			\label{eq:simp2}
			\text{tr}(M\Sigma)=\sum_i \sigma_i M_{ii}
		\end{align}
		This expression is maximized under the conditions above when $M$ is the identity.
		\begin{align}
			M=V^TR^TU=I\\
			\implies R=UV^T
		\end{align}
		To ensure that $R$ is a proper rotation, we require $|R|=1$. When $|UV^T|=-1$, we may negate the last column of $U$ to achieve this condition.
	\end{proof}
	
	\begin{theorem} (Horn algorithm)
		Given the same inputs as in Theorem \ref{thm:kabsch}, the unit quaternion corresponding to the rotation which minimizes $L(R)$ is the maximal vector of the symmetric matrix $N$, defined in Eq (\ref{eq:quatcost}).
	\end{theorem}
	\begin{proof}
		We will not repeat the entirety of Horn's derivation here \cite{Horn1987}. For a review of the quaternion algebra, refer to \S \ref{sec:quaternions}. We adapt our notation here to match the preceding sections, so that $\mathbf{v}_{A,i}$ and $\mathbf{v}_{B,i}$ are ordinary vectors in $\mathbb R^3$, and $v_{A,i}$ and $v_{B,i}$ are their quaternion counterparts using Eq (\ref{eq:quaternify}).
		\par Horn attempts to maximize the dot products of the rotated first set of vectors with the second set of vectors, i.e. the equivalent of Eq (\ref{eq:cost4}) from the proof of the Kabsch algorithm.
		\begin{align}
			f(q)=\sum_i (qv_{A,i}q^*)\cdot v_{B,i}.
		\end{align}
		In this expression, $q$ is a unit quaternion, $qv_{A,i}q^*$ is computed via the quaternion product, defined in Eq (\ref{eq:quatmul}), and $\cdot$ is the Euclidean inner product. Using the interpretation of quaternion multiplication as a matrix multiplication, shown in Eq \ref{eq:quatmatmul}), Horn switches notation, rewriting $f(q)$ as a quadratic form:
		\begin{align}
			\label{eq:quatcost}
			f(q)=\mathbf{q}^TN\mathbf{q}.
		\end{align}
		In this expression, $\mathbf{q}$ are the components of the quaternion $q$, now regarded as a column vector. $N$ is a symmetric matrix formed using the entries in $H=BA^T$ from the Kabsch derivation. The key insight now is that Eq (\ref{eq:quatcost}) can be analyzed in terms of the eigendecomposition of $N$. N is symmetric, and therefore it possesses a complete, orthonormal eigenbasis.  Since the eigenvectors ${e_i}$ are a basis for $\mathbb R^4$ the desired unit quaternion may be expressed as $\mathbf{q}=a_1\mathbf e_1+a_2\mathbf e_2+a_3\mathbf e_3+a_4\mathbf e_4$, with $a_i\in \mathbb R$. The norm $||q||_{\text{quat}}$ can be calculated as $a_1^2+a_2^2+a_3^2+a_4^2=1$, and now:
		\begin{align}
			\label{eq:horn}
			\mathbf{q}^TN\mathbf{q}=\alpha_1^2\lambda_1+\alpha_2^2\lambda_2+\alpha_3^2\lambda_3+\alpha_4^2\lambda_4,
		\end{align}
		which is maximized when the largest eigenvalue $\lambda_{\text{max}}$ and its corresponding eigenvector $\mathbf e_{\text{max}}$ are chosen as $\mathbf q$. One can therefore use a power iteration method or similar algorithm to calculate the principal eigenvalue and eigenvector of $N$ \cite{Sha2021}.
	\end{proof}
	\section{Can we use the Kabsch or Horn algorithms as inspiration for Minkowski vector alignment?}
	What follows is an informal discussion of why direct application of these algorithms is not feasible. For Minkowski vector alignment, we follow the Kabsch algorithm in defining the objective function as:
	\begin{equation}
		\label{eq:cost_lorentz}
		f(\Lambda)=\sum_i \|v_{B,i}-\Lambda v_{A,i}\|^2,
	\end{equation}
	where we again use the Frobenius norm $\|\cdot\|$, thereby aiming to match the numerical values of the 4-vector components between frames $A$ and $B$.
	\subsection{Attempt using the Kabsch algorithm}
	\begin{prop}
		A direct analogue of the Kabsch algorithm does not provide a closed-form solution to Eq (\ref{eq:cost_lorentz}) without additional theoretical and computational tools.
	\end{prop}
	The Kabsch algorithm's elegance stems from the fact that singular value decomposition of a matrix $H=BA^T$ automatically produces orthogonal matrices $U,V\in O(N)$, which allows the optimization over $SO(N)$ to be solved in closed form via Eq (\ref{eq:simp1})-Eq (\ref{eq:simp2}). For the Lorentz case, one could in principle seek to optimize over $SO(3,1)_+$ using an analogous decomposition. However, this would require a Lorentz singular value decomposition \cite{Verstraete2002a,Dieci2004}, i.e. a decomposition of $H$ into factors preserving the Minkowski metric rather than the Euclidean metric.
	\par While theoretical work on such decompositions exists, to our knowledge there is no widely available, computationally efficient implementation suitable for numerical optimization. Developing such methods is beyond the scope of this work.
	\subsection{Attempt using the Horn algorithm}
	\par Using the Horn algorithm, the situation becomes more complex (literally). Just as unit quaternions represent rotations in $SO(3)$, unit \textbf{biquaternions} represent Lorentz transformations \cite{Pellat-Finet2024}. A biquaternion is a vector space over $\mathbb C$ rather than $\mathbb R$. The rules of quaternion conjugation and multiplication remain identical. Crucially, whereas vectors in $\mathbb R^3$ were packaged into the imaginary part of the quaternion, 4-vectors must instead be packaged as Minkowski biquaternions, also known as \textbf{minquats}, so that a 4-vector $\mathbf{v}=(t,x,y,z)$ is represented as
	\begin{equation}
		v=(t,ix, iy, iz),
	\end{equation}
	where $i$ here is the imaginary unit in $\mathbb C$, so that the biquaternion $v$ has purely real (in $\mathbb C$) first component and purely imaginary (in $\mathbb C$) remaining components.
	\par Biquaternions have two types of conjugation: complex conjugation in the field, defined by
	\begin{equation}
		\bar q=(\bar q_0,\bar q_1,\bar q_2,\bar q_3),
	\end{equation}
	and conjugation over the quaternion algebra, represented by $q^*$ in Eq (\ref{eq:quatconj}).
	The two types of conjugation commute $(\bar q)^*=\overline{(q^*)}=\overline{q^*}$.
	The biquaternion norm is defined as the complex number $||q||_{\text{biquat}}=qq^*\in \mathbb C$. The action of the unit quaternion $q$ associated with Lorentz transformation $\Lambda$ on a minquat $v$ is defined by
	\begin{equation}
		\label{eq:quatboost}
		\Lambda_q(v)=qv\overline{q^*}.
	\end{equation}
	One may postulate that perhaps we can maximize a Minkowski inner product (for alignment of time-like vectors), since many of the quaternion properties still hold. Unfortunately, a unit biquaternion's components may be unbounded, for example $(i\cosh\phi,\sinh\phi,0,0)$ which has biquaternion norm $\sinh^2\phi-\cosh^2\phi=1$. The presence of negative contributions to the norm prevents us from writing down the maximum as in Eq (\ref{eq:horn}).
	\begin{prop}
		\label{thm:horn}
		Attempting to generalize the Horn approach with the biquaternion algebra does not linearize the Lorentz alignment problem into a principal eigenvector computation as it does for $SO(3)$.
	\end{prop}
	\par For $SO(3)$, the matrix $N$ is real and symmetric, the coefficients $\alpha_i$ are real, and therefore $\alpha_i^2\geq 0$. The optimization problem
	\begin{align}
		\text{maximize } \sum_i\alpha_i^2\lambda_i\quad\text{subject to}\quad\sum_i\alpha_i^2=1
	\end{align}
	is convex, since the objective is a convex combination of the eigenvalues $\lambda_i$. The maximum is achieved at a vertex of the constraint simplex, i.e., when $\alpha_j=1$ for the index $j$ corresponding to the largest eigenvalue $\lambda_{\text{max}}$, and $\alpha_i=0$ otherwise. This is precisely the principal eigenvector solution.
	
	\par For Lorentz transformations, one would work with unit biquaternions, where $N$ is Hermitian and has real eigenvalues $\lambda_i$ and orthogonal eigenvectors $e_i$. The unit biquaternion constraint is
	\begin{align}
		\label{eq:weird}
		\alpha_1^2+\alpha_2^2+\alpha_3^2+\alpha_4^2=1,
	\end{align}
	where now $\alpha_i\in\mathbb{C}$. Note that the constraint involves $\alpha_i^2$ (not $|\alpha_i|^2$) because the biquaternion norm $qq^*$ can be complex-valued; this reflects the indefinite Minkowski signature. Consequently, the optimization problem becomes
	\begin{align}
		\text{maximize } \sum_i|\alpha_i|^2\lambda_i\quad\text{subject to}\quad\sum_i\alpha_i^2=1.
	\end{align}
	Since $\alpha_i\in\mathbb{C}$, the quantities $\alpha_i^2$ can be negative or even complex (e.g., if $\alpha_i=e^{i\pi/4}$, then $\alpha_i^2=i$), and the optimization problem is no longer convex. The principal eigenvector method, which relies on the convex structure, does not apply. While one could still attempt to solve this nonlinear problem numerically, we no longer realize the computational advantage of Horn's algorithm for $SO(3)$.
	\pagebreak


\begin{thebibliography}{15}%
		\makeatletter
		\providecommand \@ifxundefined [1]{%
			\@ifx{#1\undefined}
		}%
		\providecommand \@ifnum [1]{%
			\ifnum #1\expandafter \@firstoftwo
			\else \expandafter \@secondoftwo
			\fi
		}%
		\providecommand \@ifx [1]{%
			\ifx #1\expandafter \@firstoftwo
			\else \expandafter \@secondoftwo
			\fi
		}%
		\providecommand \natexlab [1]{#1}%
		\providecommand \enquote  [1]{``#1''}%
		\providecommand \bibnamefont  [1]{#1}%
		\providecommand \bibfnamefont [1]{#1}%
		\providecommand \citenamefont [1]{#1}%
		\providecommand \href@noop [0]{\@secondoftwo}%
		\providecommand \href [0]{\begingroup \@sanitize@url \@href}%
		\providecommand \@href[1]{\@@startlink{#1}\@@href}%
		\providecommand \@@href[1]{\endgroup#1\@@endlink}%
		\providecommand \@sanitize@url [0]{\catcode `\\12\catcode `\$12\catcode
			`\&12\catcode `\#12\catcode `\^12\catcode `\_12\catcode `\%12\relax}%
		\providecommand \@@startlink[1]{}%
		\providecommand \@@endlink[0]{}%
		\providecommand \url  [0]{\begingroup\@sanitize@url \@url }%
		\providecommand \@url [1]{\endgroup\@href {#1}{\urlprefix }}%
		\providecommand \urlprefix  [0]{URL }%
		\providecommand \Eprint [0]{\href }%
		\providecommand \doibase [0]{https://doi.org/}%
		\providecommand \selectlanguage [0]{\@gobble}%
		\providecommand \bibinfo  [0]{\@secondoftwo}%
		\providecommand \bibfield  [0]{\@secondoftwo}%
		\providecommand \translation [1]{[#1]}%
		\providecommand \BibitemOpen [0]{}%
		\providecommand \bibitemStop [0]{}%
		\providecommand \bibitemNoStop [0]{.\EOS\space}%
		\providecommand \EOS [0]{\spacefactor3000\relax}%
		\providecommand \BibitemShut  [1]{\csname bibitem#1\endcsname}%
		\let\auto@bib@innerbib\@empty
		\bibitem [{\citenamefont {Kabsch}(1976)}]{Kabsch1976}%
		\BibitemOpen
		\bibfield  {author} {\bibinfo {author} {\bibfnamefont {W.}~\bibnamefont
				{Kabsch}},\ }\bibfield  {title} {\bibinfo {title} {{A solution for the best
					rotation to relate two sets of vectors}},\ }\bibfield  {journal} {\bibinfo
			{journal} {Acta Crystallographica Section A}\ }\textbf {\bibinfo {volume}
			{32}},\ \href {https://doi.org/10.1107/S0567739476001873}
		{10.1107/S0567739476001873} (\bibinfo {year} {1976})\BibitemShut {NoStop}%
		\bibitem [{\citenamefont {Horn}(1987)}]{Horn1987}%
		\BibitemOpen
		\bibfield  {author} {\bibinfo {author} {\bibfnamefont {B.~K.~P.}\
				\bibnamefont {Horn}},\ }\bibfield  {title} {\bibinfo {title} {{Closed-form
					solution of absolute orientation using unit quaternions}},\ }\bibfield
		{journal} {\bibinfo  {journal} {Journal of the Optical Society of America A}\
		}\textbf {\bibinfo {volume} {4}},\ \href
		{https://doi.org/10.1364/josaa.4.000629} {10.1364/josaa.4.000629} (\bibinfo
		{year} {1987})\BibitemShut {NoStop}%
		\bibitem [{\citenamefont {Hall}(2015)}]{Hall2013a}%
		\BibitemOpen
		\bibfield  {author} {\bibinfo {author} {\bibfnamefont {B.~C.}\ \bibnamefont
				{Hall}},\ }\href {https://doi.org/10.1007/978-1-4614-7116-5_16} {\emph
			{\bibinfo {title} {{Lie Groups, Lie Algebras, and Representations: An
						Elementary Introduction}}}},\ \bibinfo {edition} {2nd}\ ed.\ (\bibinfo
		{publisher} {Springer Cham},\ \bibinfo {address} {Cham, Switzerland},\
		\bibinfo {year} {2015})\BibitemShut {NoStop}%
		\bibitem [{\citenamefont {Dai}(2015)}]{Dai2015}%
		\BibitemOpen
		\bibfield  {author} {\bibinfo {author} {\bibfnamefont {J.~S.}\ \bibnamefont
				{Dai}},\ }\bibfield  {title} {\bibinfo {title} {{Euler–Rodrigues formula
					variations, quaternion conjugation and intrinsic connections}},\ }\href
		{https://doi.org/10.1016/j.mechmachtheory.2015.03.004} {\bibfield  {journal}
			{\bibinfo  {journal} {Mechanism and Machine Theory}\ }\textbf {\bibinfo
				{volume} {92}},\ \bibinfo {pages} {144} (\bibinfo {year} {2015})}\BibitemShut
		{NoStop}%
		\bibitem [{\citenamefont {Haber}(2024)}]{Haber2024}%
		\BibitemOpen
		\bibfield  {author} {\bibinfo {author} {\bibfnamefont {H.~E.}\ \bibnamefont
				{Haber}},\ }\bibfield  {title} {\bibinfo {title} {{Explicit Form for the Most
					General Lorentz Transformation Revisited}},\ }\href
		{https://doi.org/10.3390/SYM16091155} {\bibfield  {journal} {\bibinfo
				{journal} {Symmetry 2024, Vol. 16, Page 1155}\ }\textbf {\bibinfo {volume}
				{16}},\ \bibinfo {pages} {1155} (\bibinfo {year} {2024})},\ \Eprint
		{https://arxiv.org/abs/2312.12969} {arXiv:2312.12969} \BibitemShut {NoStop}%
		\bibitem [{\citenamefont {Burden}\ \emph {et~al.}(2015)\citenamefont {Burden},
			\citenamefont {{L. Burden}},\ and\ \citenamefont {{Douglas
					Faires}}}]{Burden2016}%
		\BibitemOpen
		\bibfield  {author} {\bibinfo {author} {\bibfnamefont {A.}~\bibnamefont
				{Burden}}, \bibinfo {author} {\bibfnamefont {R.}~\bibnamefont {{L.
						Burden}}},\ and\ \bibinfo {author} {\bibfnamefont {J.}~\bibnamefont {{Douglas
						Faires}}},\ }\href@noop {} {\emph {\bibinfo {title} {Cengage Learning}}}\
		(\bibinfo {year} {2015})\BibitemShut {NoStop}%
		\bibitem [{\citenamefont {Virtanen}\ \emph {et~al.}(2020)\citenamefont
			{Virtanen}, \citenamefont {Gommers}, \citenamefont {Oliphant}, \citenamefont
			{Haberland}, \citenamefont {Reddy}, \citenamefont {Cournapeau}, \citenamefont
			{Burovski}, \citenamefont {Peterson}, \citenamefont {Weckesser},
			\citenamefont {Bright}, \citenamefont {van~der Walt}, \citenamefont {Brett},
			\citenamefont {Wilson}, \citenamefont {Millman}, \citenamefont {Mayorov},
			\citenamefont {Nelson}, \citenamefont {Jones}, \citenamefont {Kern},
			\citenamefont {Larson}, \citenamefont {Carey}, \citenamefont {Polat},
			\citenamefont {Feng}, \citenamefont {Moore}, \citenamefont {VanderPlas},
			\citenamefont {Laxalde}, \citenamefont {Perktold}, \citenamefont {Cimrman},
			\citenamefont {Henriksen}, \citenamefont {Quintero}, \citenamefont {Harris},
			\citenamefont {Archibald}, \citenamefont {Ribeiro}, \citenamefont
			{Pedregosa}, \citenamefont {van Mulbregt}, \citenamefont {Vijaykumar},
			\citenamefont {Bardelli}, \citenamefont {Rothberg}, \citenamefont {Hilboll},
			\citenamefont {Kloeckner}, \citenamefont {Scopatz}, \citenamefont {Lee},
			\citenamefont {Rokem}, \citenamefont {Woods}, \citenamefont {Fulton},
			\citenamefont {Masson}, \citenamefont {H{\"{a}}ggstr{\"{o}}m}, \citenamefont
			{Fitzgerald}, \citenamefont {Nicholson}, \citenamefont {Hagen}, \citenamefont
			{Pasechnik}, \citenamefont {Olivetti}, \citenamefont {Martin}, \citenamefont
			{Wieser}, \citenamefont {Silva}, \citenamefont {Lenders}, \citenamefont
			{Wilhelm}, \citenamefont {Young}, \citenamefont {Price}, \citenamefont
			{Ingold}, \citenamefont {Allen}, \citenamefont {Lee}, \citenamefont {Audren},
			\citenamefont {Probst}, \citenamefont {Dietrich}, \citenamefont {Silterra},
			\citenamefont {Webber}, \citenamefont {Slavi{\v{c}}}, \citenamefont
			{Nothman}, \citenamefont {Buchner}, \citenamefont {Kulick}, \citenamefont
			{Sch{\"{o}}nberger}, \citenamefont {{de Miranda Cardoso}}, \citenamefont
			{Reimer}, \citenamefont {Harrington}, \citenamefont {Rodr{\'{i}}guez},
			\citenamefont {Nunez-Iglesias}, \citenamefont {Kuczynski}, \citenamefont
			{Tritz}, \citenamefont {Thoma}, \citenamefont {Newville}, \citenamefont
			{K{\"{u}}mmerer}, \citenamefont {Bolingbroke}, \citenamefont {Tartre},
			\citenamefont {Pak}, \citenamefont {Smith}, \citenamefont {Nowaczyk},
			\citenamefont {Shebanov}, \citenamefont {Pavlyk}, \citenamefont {Brodtkorb},
			\citenamefont {Lee}, \citenamefont {McGibbon}, \citenamefont {Feldbauer},
			\citenamefont {Lewis}, \citenamefont {Tygier}, \citenamefont {Sievert},
			\citenamefont {Vigna}, \citenamefont {Peterson}, \citenamefont {More},
			\citenamefont {Pudlik}, \citenamefont {Oshima}, \citenamefont {Pingel},
			\citenamefont {Robitaille}, \citenamefont {Spura}, \citenamefont {Jones},
			\citenamefont {Cera}, \citenamefont {Leslie}, \citenamefont {Zito},
			\citenamefont {Krauss}, \citenamefont {Upadhyay}, \citenamefont {Halchenko},\
			and\ \citenamefont {V{\'{a}}zquez-Baeza}}]{Virtanen2020a}%
		\BibitemOpen
		\bibfield  {author} {\bibinfo {author} {\bibfnamefont {P.}~\bibnamefont
				{Virtanen}}, \bibinfo {author} {\bibfnamefont {R.}~\bibnamefont {Gommers}},
			\bibinfo {author} {\bibfnamefont {T.~E.}\ \bibnamefont {Oliphant}}, \bibinfo
			{author} {\bibfnamefont {M.}~\bibnamefont {Haberland}}, \bibinfo {author}
			{\bibfnamefont {T.}~\bibnamefont {Reddy}}, \bibinfo {author} {\bibfnamefont
				{D.}~\bibnamefont {Cournapeau}}, \bibinfo {author} {\bibfnamefont
				{E.}~\bibnamefont {Burovski}}, \bibinfo {author} {\bibfnamefont
				{P.}~\bibnamefont {Peterson}}, \bibinfo {author} {\bibfnamefont
				{W.}~\bibnamefont {Weckesser}}, \bibinfo {author} {\bibfnamefont
				{J.}~\bibnamefont {Bright}}, \bibinfo {author} {\bibfnamefont {S.~J.}\
				\bibnamefont {van~der Walt}}, \bibinfo {author} {\bibfnamefont
				{M.}~\bibnamefont {Brett}}, \bibinfo {author} {\bibfnamefont
				{J.}~\bibnamefont {Wilson}}, \bibinfo {author} {\bibfnamefont {K.~J.}\
				\bibnamefont {Millman}}, \bibinfo {author} {\bibfnamefont {N.}~\bibnamefont
				{Mayorov}}, \bibinfo {author} {\bibfnamefont {A.~R.~J.}\ \bibnamefont
				{Nelson}}, \bibinfo {author} {\bibfnamefont {E.}~\bibnamefont {Jones}},
			\bibinfo {author} {\bibfnamefont {R.}~\bibnamefont {Kern}}, \bibinfo {author}
			{\bibfnamefont {E.}~\bibnamefont {Larson}}, \bibinfo {author} {\bibfnamefont
				{C.~J.}\ \bibnamefont {Carey}}, \bibinfo {author} {\bibfnamefont
				{I.}~\bibnamefont {Polat}}, \bibinfo {author} {\bibfnamefont
				{Y.}~\bibnamefont {Feng}}, \bibinfo {author} {\bibfnamefont {E.~W.}\
				\bibnamefont {Moore}}, \bibinfo {author} {\bibfnamefont {J.}~\bibnamefont
				{VanderPlas}}, \bibinfo {author} {\bibfnamefont {D.}~\bibnamefont {Laxalde}},
			\bibinfo {author} {\bibfnamefont {J.}~\bibnamefont {Perktold}}, \bibinfo
			{author} {\bibfnamefont {R.}~\bibnamefont {Cimrman}}, \bibinfo {author}
			{\bibfnamefont {I.}~\bibnamefont {Henriksen}}, \bibinfo {author}
			{\bibfnamefont {E.~A.}\ \bibnamefont {Quintero}}, \bibinfo {author}
			{\bibfnamefont {C.~R.}\ \bibnamefont {Harris}}, \bibinfo {author}
			{\bibfnamefont {A.~M.}\ \bibnamefont {Archibald}}, \bibinfo {author}
			{\bibfnamefont {A.~H.}\ \bibnamefont {Ribeiro}}, \bibinfo {author}
			{\bibfnamefont {F.}~\bibnamefont {Pedregosa}}, \bibinfo {author}
			{\bibfnamefont {P.}~\bibnamefont {van Mulbregt}}, \bibinfo {author}
			{\bibfnamefont {A.}~\bibnamefont {Vijaykumar}}, \bibinfo {author}
			{\bibfnamefont {A.~P.}\ \bibnamefont {Bardelli}}, \bibinfo {author}
			{\bibfnamefont {A.}~\bibnamefont {Rothberg}}, \bibinfo {author}
			{\bibfnamefont {A.}~\bibnamefont {Hilboll}}, \bibinfo {author} {\bibfnamefont
				{A.}~\bibnamefont {Kloeckner}}, \bibinfo {author} {\bibfnamefont
				{A.}~\bibnamefont {Scopatz}}, \bibinfo {author} {\bibfnamefont
				{A.}~\bibnamefont {Lee}}, \bibinfo {author} {\bibfnamefont {A.}~\bibnamefont
				{Rokem}}, \bibinfo {author} {\bibfnamefont {C.~N.}\ \bibnamefont {Woods}},
			\bibinfo {author} {\bibfnamefont {C.}~\bibnamefont {Fulton}}, \bibinfo
			{author} {\bibfnamefont {C.}~\bibnamefont {Masson}}, \bibinfo {author}
			{\bibfnamefont {C.}~\bibnamefont {H{\"{a}}ggstr{\"{o}}m}}, \bibinfo {author}
			{\bibfnamefont {C.}~\bibnamefont {Fitzgerald}}, \bibinfo {author}
			{\bibfnamefont {D.~A.}\ \bibnamefont {Nicholson}}, \bibinfo {author}
			{\bibfnamefont {D.~R.}\ \bibnamefont {Hagen}}, \bibinfo {author}
			{\bibfnamefont {D.~V.}\ \bibnamefont {Pasechnik}}, \bibinfo {author}
			{\bibfnamefont {E.}~\bibnamefont {Olivetti}}, \bibinfo {author}
			{\bibfnamefont {E.}~\bibnamefont {Martin}}, \bibinfo {author} {\bibfnamefont
				{E.}~\bibnamefont {Wieser}}, \bibinfo {author} {\bibfnamefont
				{F.}~\bibnamefont {Silva}}, \bibinfo {author} {\bibfnamefont
				{F.}~\bibnamefont {Lenders}}, \bibinfo {author} {\bibfnamefont
				{F.}~\bibnamefont {Wilhelm}}, \bibinfo {author} {\bibfnamefont
				{G.}~\bibnamefont {Young}}, \bibinfo {author} {\bibfnamefont {G.~A.}\
				\bibnamefont {Price}}, \bibinfo {author} {\bibfnamefont {G.-L.}\ \bibnamefont
				{Ingold}}, \bibinfo {author} {\bibfnamefont {G.~E.}\ \bibnamefont {Allen}},
			\bibinfo {author} {\bibfnamefont {G.~R.}\ \bibnamefont {Lee}}, \bibinfo
			{author} {\bibfnamefont {H.}~\bibnamefont {Audren}}, \bibinfo {author}
			{\bibfnamefont {I.}~\bibnamefont {Probst}}, \bibinfo {author} {\bibfnamefont
				{J.~P.}\ \bibnamefont {Dietrich}}, \bibinfo {author} {\bibfnamefont
				{J.}~\bibnamefont {Silterra}}, \bibinfo {author} {\bibfnamefont {J.~T.}\
				\bibnamefont {Webber}}, \bibinfo {author} {\bibfnamefont {J.}~\bibnamefont
				{Slavi{\v{c}}}}, \bibinfo {author} {\bibfnamefont {J.}~\bibnamefont
				{Nothman}}, \bibinfo {author} {\bibfnamefont {J.}~\bibnamefont {Buchner}},
			\bibinfo {author} {\bibfnamefont {J.}~\bibnamefont {Kulick}}, \bibinfo
			{author} {\bibfnamefont {J.~L.}\ \bibnamefont {Sch{\"{o}}nberger}}, \bibinfo
			{author} {\bibfnamefont {J.~V.}\ \bibnamefont {{de Miranda Cardoso}}},
			\bibinfo {author} {\bibfnamefont {J.}~\bibnamefont {Reimer}}, \bibinfo
			{author} {\bibfnamefont {J.}~\bibnamefont {Harrington}}, \bibinfo {author}
			{\bibfnamefont {J.~L.~C.}\ \bibnamefont {Rodr{\'{i}}guez}}, \bibinfo {author}
			{\bibfnamefont {J.}~\bibnamefont {Nunez-Iglesias}}, \bibinfo {author}
			{\bibfnamefont {J.}~\bibnamefont {Kuczynski}}, \bibinfo {author}
			{\bibfnamefont {K.}~\bibnamefont {Tritz}}, \bibinfo {author} {\bibfnamefont
				{M.}~\bibnamefont {Thoma}}, \bibinfo {author} {\bibfnamefont
				{M.}~\bibnamefont {Newville}}, \bibinfo {author} {\bibfnamefont
				{M.}~\bibnamefont {K{\"{u}}mmerer}}, \bibinfo {author} {\bibfnamefont
				{M.}~\bibnamefont {Bolingbroke}}, \bibinfo {author} {\bibfnamefont
				{M.}~\bibnamefont {Tartre}}, \bibinfo {author} {\bibfnamefont
				{M.}~\bibnamefont {Pak}}, \bibinfo {author} {\bibfnamefont {N.~J.}\
				\bibnamefont {Smith}}, \bibinfo {author} {\bibfnamefont {N.}~\bibnamefont
				{Nowaczyk}}, \bibinfo {author} {\bibfnamefont {N.}~\bibnamefont {Shebanov}},
			\bibinfo {author} {\bibfnamefont {O.}~\bibnamefont {Pavlyk}}, \bibinfo
			{author} {\bibfnamefont {P.~A.}\ \bibnamefont {Brodtkorb}}, \bibinfo {author}
			{\bibfnamefont {P.}~\bibnamefont {Lee}}, \bibinfo {author} {\bibfnamefont
				{R.~T.}\ \bibnamefont {McGibbon}}, \bibinfo {author} {\bibfnamefont
				{R.}~\bibnamefont {Feldbauer}}, \bibinfo {author} {\bibfnamefont
				{S.}~\bibnamefont {Lewis}}, \bibinfo {author} {\bibfnamefont
				{S.}~\bibnamefont {Tygier}}, \bibinfo {author} {\bibfnamefont
				{S.}~\bibnamefont {Sievert}}, \bibinfo {author} {\bibfnamefont
				{S.}~\bibnamefont {Vigna}}, \bibinfo {author} {\bibfnamefont
				{S.}~\bibnamefont {Peterson}}, \bibinfo {author} {\bibfnamefont
				{S.}~\bibnamefont {More}}, \bibinfo {author} {\bibfnamefont {T.}~\bibnamefont
				{Pudlik}}, \bibinfo {author} {\bibfnamefont {T.}~\bibnamefont {Oshima}},
			\bibinfo {author} {\bibfnamefont {T.~J.}\ \bibnamefont {Pingel}}, \bibinfo
			{author} {\bibfnamefont {T.~P.}\ \bibnamefont {Robitaille}}, \bibinfo
			{author} {\bibfnamefont {T.}~\bibnamefont {Spura}}, \bibinfo {author}
			{\bibfnamefont {T.~R.}\ \bibnamefont {Jones}}, \bibinfo {author}
			{\bibfnamefont {T.}~\bibnamefont {Cera}}, \bibinfo {author} {\bibfnamefont
				{T.}~\bibnamefont {Leslie}}, \bibinfo {author} {\bibfnamefont
				{T.}~\bibnamefont {Zito}}, \bibinfo {author} {\bibfnamefont {T.}~\bibnamefont
				{Krauss}}, \bibinfo {author} {\bibfnamefont {U.}~\bibnamefont {Upadhyay}},
			\bibinfo {author} {\bibfnamefont {Y.~O.}\ \bibnamefont {Halchenko}},\ and\
			\bibinfo {author} {\bibfnamefont {Y.}~\bibnamefont {V{\'{a}}zquez-Baeza}},\
		}\bibfield  {title} {\bibinfo {title} {{SciPy 1.0: fundamental algorithms for
					scientific computing in Python}},\ }\href
		{https://doi.org/10.1038/s41592-019-0686-2} {\bibfield  {journal} {\bibinfo
				{journal} {Nature Methods}\ }\textbf {\bibinfo {volume} {17}},\ \bibinfo
			{pages} {261} (\bibinfo {year} {2020})}\BibitemShut {NoStop}%
		\bibitem [{\citenamefont {Strang}(2003)}]{Strang2003}%
		\BibitemOpen
		\bibfield  {author} {\bibinfo {author} {\bibfnamefont {G.}~\bibnamefont
				{Strang}},\ }\bibfield  {title} {\bibinfo {title} {{Linear Algebra and Its
					Applications}},\ }\bibfield  {journal} {\bibinfo  {journal} {Linear Algebra
				and Its Applications}\ }\textbf {\bibinfo {volume} {373}},\ \href
		{https://doi.org/10.2307/2978065} {10.2307/2978065} (\bibinfo {year}
		{2003})\BibitemShut {NoStop}%
		\bibitem [{\citenamefont {Harris}\ \emph {et~al.}(2020)\citenamefont {Harris},
			\citenamefont {Millman}, \citenamefont {van~der Walt}, \citenamefont
			{Gommers}, \citenamefont {Virtanen}, \citenamefont {Cournapeau},
			\citenamefont {Wieser}, \citenamefont {Taylor}, \citenamefont {Berg},
			\citenamefont {Smith}, \citenamefont {Kern}, \citenamefont {Picus},
			\citenamefont {Hoyer}, \citenamefont {van Kerkwijk}, \citenamefont {Brett},
			\citenamefont {Haldane}, \citenamefont {del R{\'{i}}o}, \citenamefont
			{Wiebe}, \citenamefont {Peterson}, \citenamefont {G{\'{e}}rard-Marchant},
			\citenamefont {Sheppard}, \citenamefont {Reddy}, \citenamefont {Weckesser},
			\citenamefont {Abbasi}, \citenamefont {Gohlke},\ and\ \citenamefont
			{Oliphant}}]{Harris2020}%
		\BibitemOpen
		\bibfield  {author} {\bibinfo {author} {\bibfnamefont {C.~R.}\ \bibnamefont
				{Harris}}, \bibinfo {author} {\bibfnamefont {K.~J.}\ \bibnamefont {Millman}},
			\bibinfo {author} {\bibfnamefont {S.~J.}\ \bibnamefont {van~der Walt}},
			\bibinfo {author} {\bibfnamefont {R.}~\bibnamefont {Gommers}}, \bibinfo
			{author} {\bibfnamefont {P.}~\bibnamefont {Virtanen}}, \bibinfo {author}
			{\bibfnamefont {D.}~\bibnamefont {Cournapeau}}, \bibinfo {author}
			{\bibfnamefont {E.}~\bibnamefont {Wieser}}, \bibinfo {author} {\bibfnamefont
				{J.}~\bibnamefont {Taylor}}, \bibinfo {author} {\bibfnamefont
				{S.}~\bibnamefont {Berg}}, \bibinfo {author} {\bibfnamefont {N.~J.}\
				\bibnamefont {Smith}}, \bibinfo {author} {\bibfnamefont {R.}~\bibnamefont
				{Kern}}, \bibinfo {author} {\bibfnamefont {M.}~\bibnamefont {Picus}},
			\bibinfo {author} {\bibfnamefont {S.}~\bibnamefont {Hoyer}}, \bibinfo
			{author} {\bibfnamefont {M.~H.}\ \bibnamefont {van Kerkwijk}}, \bibinfo
			{author} {\bibfnamefont {M.}~\bibnamefont {Brett}}, \bibinfo {author}
			{\bibfnamefont {A.}~\bibnamefont {Haldane}}, \bibinfo {author} {\bibfnamefont
				{J.~F.}\ \bibnamefont {del R{\'{i}}o}}, \bibinfo {author} {\bibfnamefont
				{M.}~\bibnamefont {Wiebe}}, \bibinfo {author} {\bibfnamefont
				{P.}~\bibnamefont {Peterson}}, \bibinfo {author} {\bibfnamefont
				{P.}~\bibnamefont {G{\'{e}}rard-Marchant}}, \bibinfo {author} {\bibfnamefont
				{K.}~\bibnamefont {Sheppard}}, \bibinfo {author} {\bibfnamefont
				{T.}~\bibnamefont {Reddy}}, \bibinfo {author} {\bibfnamefont
				{W.}~\bibnamefont {Weckesser}}, \bibinfo {author} {\bibfnamefont
				{H.}~\bibnamefont {Abbasi}}, \bibinfo {author} {\bibfnamefont
				{C.}~\bibnamefont {Gohlke}},\ and\ \bibinfo {author} {\bibfnamefont {T.~E.}\
				\bibnamefont {Oliphant}},\ }\href {https://doi.org/10.1038/s41586-020-2649-2}
		{\bibinfo {title} {{Array programming with NumPy}}} (\bibinfo {year}
		{2020}),\ \Eprint {https://arxiv.org/abs/2006.10256} {arXiv:2006.10256}
		\BibitemShut {NoStop}%
		\bibitem [{\citenamefont {Brewin}(2011)}]{Brewin2011}%
		\BibitemOpen
		\bibfield  {author} {\bibinfo {author} {\bibfnamefont {L.}~\bibnamefont
				{Brewin}},\ }\bibfield  {title} {\bibinfo {title} {{An Einstein-Bianchi
					system for Smooth Lattice General Relativity. I. The Schwarzschild
					spacetime}},\ }\bibfield  {journal} {\bibinfo  {journal} {Physical Review D -
				Particles, Fields, Gravitation and Cosmology}\ }\textbf {\bibinfo {volume}
			{85}},\ \href {https://doi.org/10.1103/PhysRevD.85.124045}
		{10.1103/PhysRevD.85.124045} (\bibinfo {year} {2011}),\ \Eprint
		{https://arxiv.org/abs/1101.3171v2} {arXiv:1101.3171v2} \BibitemShut
		{NoStop}%
		\bibitem [{\citenamefont {Al-Mohy}\ and\ \citenamefont
			{Higham}(2012)}]{Al-Mohy2012}%
		\BibitemOpen
		\bibfield  {author} {\bibinfo {author} {\bibfnamefont {A.~H.}\ \bibnamefont
				{Al-Mohy}}\ and\ \bibinfo {author} {\bibfnamefont {N.~J.}\ \bibnamefont
				{Higham}},\ }\bibfield  {title} {\bibinfo {title} {{Improved inverse scaling
					and Squaring algorithms for the matrix logarithm}},\ }\bibfield  {journal}
		{\bibinfo  {journal} {SIAM Journal on Scientific Computing}\ }\textbf
		{\bibinfo {volume} {34}},\ \href {https://doi.org/10.1137/110852553}
		{10.1137/110852553} (\bibinfo {year} {2012})\BibitemShut {NoStop}%
		\bibitem [{\citenamefont {Sha}\ and\ \citenamefont
			{Dokholyan}(2021)}]{Sha2021}%
		\BibitemOpen
		\bibfield  {author} {\bibinfo {author} {\bibfnamefont {C.~M.}\ \bibnamefont
				{Sha}}\ and\ \bibinfo {author} {\bibfnamefont {N.~V.}\ \bibnamefont
				{Dokholyan}},\ }\bibfield  {title} {\bibinfo {title} {{Simple exponential
					acceleration of the power iteration algorithm}}\ }\href
		{https://doi.org/10.48550/arxiv.2109.10884} {10.48550/arxiv.2109.10884}
		(\bibinfo {year} {2021}),\ \Eprint {https://arxiv.org/abs/2109.10884}
		{arXiv:2109.10884} \BibitemShut {NoStop}%
		\bibitem [{\citenamefont {Verstraete}\ \emph {et~al.}(2002)\citenamefont
			{Verstraete}, \citenamefont {Dehaene},\ and\ \citenamefont {{De
					Moor}}}]{Verstraete2002a}%
		\BibitemOpen
		\bibfield  {author} {\bibinfo {author} {\bibfnamefont {F.}~\bibnamefont
				{Verstraete}}, \bibinfo {author} {\bibfnamefont {J.}~\bibnamefont
				{Dehaene}},\ and\ \bibinfo {author} {\bibfnamefont {B.}~\bibnamefont {{De
						Moor}}},\ }\bibfield  {title} {\bibinfo {title} {{Lorentz singular-value
					decomposition and its applications to pure states of three qubits}},\
		}\bibfield  {journal} {\bibinfo  {journal} {Physical Review A - Atomic,
				Molecular, and Optical Physics}\ }\textbf {\bibinfo {volume} {65}},\ \href
		{https://doi.org/10.1103/PhysRevA.65.032308} {10.1103/PhysRevA.65.032308}
		(\bibinfo {year} {2002})\BibitemShut {NoStop}%
		\bibitem [{\citenamefont {Dieci}\ \emph {et~al.}(2004)\citenamefont {Dieci},
			\citenamefont {Elia},\ and\ \citenamefont {Lopez}}]{Dieci2004}%
		\BibitemOpen
		\bibfield  {author} {\bibinfo {author} {\bibfnamefont {L.}~\bibnamefont
				{Dieci}}, \bibinfo {author} {\bibfnamefont {C.}~\bibnamefont {Elia}},\ and\
			\bibinfo {author} {\bibfnamefont {L.}~\bibnamefont {Lopez}},\ }\bibfield
		{title} {\bibinfo {title} {{Smooth SVD on the Lorentz group with application
					to computation of Lyapunov exponents}},\ }\bibfield  {journal} {\bibinfo
			{journal} {Journal of Computational and Applied Mathematics}\ }\textbf
		{\bibinfo {volume} {164-165}},\ \href
		{https://doi.org/10.1016/S0377-0427(03)00644-7}
		{10.1016/S0377-0427(03)00644-7} (\bibinfo {year} {2004})\BibitemShut
		{NoStop}%
		\bibitem [{\citenamefont {Pellat-Finet}(2024)}]{Pellat-Finet2024}%
		\BibitemOpen
		\bibfield  {author} {\bibinfo {author} {\bibfnamefont {P.}~\bibnamefont
				{Pellat-Finet}},\ }\bibfield  {title} {\bibinfo {title} {{Quaternions and
					rotations: applications to Minkowski's four-vectors, electromagnetic waves
					and polarization optics}},\ }\href {https://arxiv.org/pdf/2407.11953} {\
			(\bibinfo {year} {2024})},\ \Eprint {https://arxiv.org/abs/2407.11953}
		{arXiv:2407.11953} \BibitemShut {NoStop}%
	\end{thebibliography}
	\end{document}